\documentclass[preprint,12pt]{elsarticle}
\usepackage{amsmath,amsfonts,amsthm}

\usepackage{graphicx}
\usepackage{float}
\usepackage{hyperref}
\usepackage[titletoc]{appendix}

\newcommand{\pa}[2]{\ensuremath{\frac{\partial #1}{\partial #2}}}

\usepackage{todonotes}
\usepackage{xcolor}

\newcommand{\mydrafttext}{}
\newcommand{\drafttext}[1]{\renewcommand{\mydrafttext}{#1}}
\newboolean{draft}  
\setboolean{draft}{true}
\ifthenelse{\boolean{draft}}
{
    \newcounter{comments}
    \drafttext{{\color{red}This document is a draft.}}
    \newcommand{\shane}[1]{\addtocounter{comments}{1}{\color{red}\bf [Shane comment \thecomments: #1]}}
    \newcommand{\erik}[1]{\addtocounter{comments}{1}{\color{blue}\bf [Erik comment \thecomments: #1]}}
}
{
\newcommand{\shane}[1]{}
\newcommand{\erik}[1]{}
}

\begin{document}

\begin{frontmatter}
	\title{Is the Finite-Time Lyapunov Exponent Field a Koopman Eigenfunction?}

	\author{Erik M. Bollt}
	\address{Electrical \& Computer Engineering \& $C^3 S^2$, the Clarkson Center for Complex Systems Science, Clarkson University \\ Potsdam, New York, 13699, USA} 
	\ead{ebollt@clarkson.edu}
	\author{Shane D. Ross}
	\address{Aerospace and Ocean Engineering, Virginia Tech \\ Blacksburg, Virginia, 24061, USA}
	\ead{sdross@vt.edu}
\begin{abstract}
This work serves as a bridge between two approaches to analysis of dynamical systems: the local, geometric analysis and the global, operator
theoretic, Koopman analysis. 
We explicitly construct vector fields where the instantaneous Lyapunov exponent field is a Koopman eigenfunction.
Restricting ourselves to polynomial vector fields to make this construction easier, we find that such vector fields do exist, and we 
explore whether such vector fields have a special structure, thus making a link between the geometric theory and the transfer operator theory.
\end{abstract}
	
\begin{keyword}
Koopman operator, spectral analysis, invariant manifolds, lyapunov exponent, dynamical systems
\end{keyword}
\date{This version: \today}

\end{frontmatter}

\section{Significance}


Two approaches to analyzing dynamical systems are the geometric approach and operator theoretic approach, exemplified in recent years by invariant manifolds and Koopman operators.
The geometric invariant manifold approach is closely related to the instantaneous version of the finite-time Lyapunov exponent (FTLE) field.
The very different, spectral and measure-based operator theoretic  approach of evolution operators, a.k.a. ``Koopmanism" involves 
 Koopman eigenfunctions (KEIGs).

In this paper, we ask a simple question, ``Is the FTLE field a KEIG?''
The answer is:\ in general, no.
This motivates the explicit construction of vector fields where the answer is yes, in the sense that the FTLE field in the infinitesimal time limit, i.e., the instantaneous Lyapunov exponent (iLE) field, is a KEIG.
Restricting ourselves to polynomial vector fields to make this construction easier, we indeed find that such vector fields do exist, 
and we explore whether such vector fields have a special structure.

\section{Instantaneous Lyapunov Exponent Analysis}\label{s:iLE}

Solutions of non-autonomous vector fields, 
\begin{equation}
     \mathbf{\dot x} = \mathbf{v}(\mathbf{x},t), \quad \mathbf{x}\in M\subset\mathbb{R}^n
    \label{vield}
\end{equation}
such as the motion of a fluid element $\mathbf{x}(t)$ in the time-dependent fluid velocity field $\mathbf{v}(\mathbf{x},t)$, can be challenging to analyze. 
The method of Lagrangian coherent structures (LCS) has 
become a popular tool to analyze structures in the phase space of low-dimensional non-autonomous vector fields \cite{HaYu2000, shadden2005definition, Haller2015}. 
A common computational framework to obtain LCS is to consider a scalar field derived from numerically integrated  trajectories (i.e., numerical approximation of solutions of \eqref{vield}, called the ``Lagrangian'' point of view in the fluid literature), in particular  the finite-time Lyapunov exponent (FTLE) field \cite{shadden2005definition,lekien2010computation}.

New  tools have recently been developed which use vector field gradients, instead of integrating  trajectories, the instantaneous limit of FTLE and LCS.  
These vector gradients are assembled into the Eulerian rate-of-strain tensor, 
\begin{equation}
    \mathbf{S}(\mathbf{x},t) = \tfrac{1}{2} \left( \nabla\mathbf{v}(\mathbf{x},t) + \nabla\mathbf{v}(\mathbf{x},t)^{T} \right)\label{eq:EST}.
\end{equation}
so-called because of its relationship to the space-fixed, or ``Eulerian'' point-of-view in the fluid literature.
The Eulerian rate-of-strain tensor was shown to provide an instantaneous approximation of LCS in two-dimensional fluid flows \cite{serra2016objective}.

Further work on this topic extended the ideas to $n$ dimensions \cite{nolan2020finite},   
showing that the minimum and maximum eigenvalues of the Eulerian rate-of-strain tensor, $s_1$ and $s_n$, 
are the limits of the backward-time and forward-time FTLE fields, 
respectively, as the  integration time goes to zero. 
Trenches of the minimum eigenvalue field can be identified as instantaneous attracting LCSs whereas ridges of the maximum eigenvalue field can be identified as instantaneous repelling LCSs. 
For the remainder of the paper, we shall refer to the minimum and maximum eigenvalues of the $\mathbf{S}$ 
as the attraction and repulsion rates, respectively, and both can be considered as 
{\it instantaneous}, 
rather than finite-time, Lyapunov exponents. 
For brevity we refer to the iLE field, the instantaneous Lyapunov exponent, field, following \cite{nolan2020finite}.

\begin{figure}[!t]
    \centering
    \includegraphics[width=0.7\linewidth]{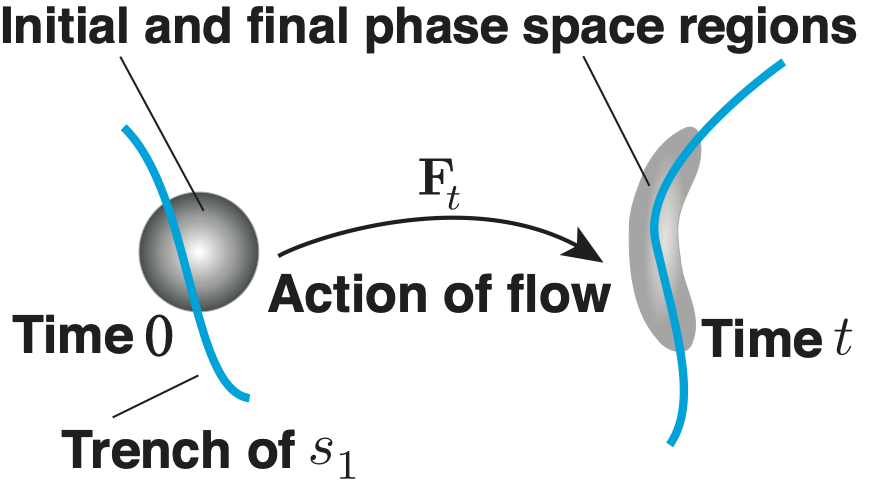}
    \caption{Schematic of the effect of the co-dimension 1 manifold, corresponding to the trench of  the attraction rate, $s_{1}$, at an initial time,  on a small parcel of phase space. The image of the trench under the time-$t$ flow map, $\mathbf{F}_t$, is also shown.
     }
    \label{fig:LCS}
\end{figure}

Broadly speaking, the trenches of the attraction rate field reveal where phase space regions will congregate under the flow, as shown in the schematic of Figure \ref{fig:LCS}.
These are instantaneous attracting LCSs, that is, the instantaneously most attracting co-dimension 1 surfaces in $M$, also referred to as instantaneous Lyapunov exponent structures (iLES) \cite{nolan2020finite}. 
In practical applications, the primary focus has been on the attraction rate given its importance for prediction \cite{nolan2018coordinated}
as opposed to repelling features which particles will diverge from before flowing independent of those features.

In the remainder of this paper, we will restrict ourselves to autonomous systems, where the vector field \eqref{vield} is independent of time, that is,
\begin{equation}
     \mathbf{\dot x} = \mathbf{v}(\mathbf{x}), \quad \mathbf{x}\in M\subset\mathbb{R}^n
    \label{vield_auton}
\end{equation}

\subsection{Instantaneous Attraction and Repulsion Rates}

For ease of exposition, we  limit our discussion in this section to
autonomous two-dimensional vector fields.
To calculate the attraction rate we first needed to calculate the gradient of the  vector field, $\nabla\mathbf{v}(\mathbf{x})$, for the   vector field $\mathbf{v}(\mathbf{x})=(u(x,y),v(x,y))$, where $u$ is the first component and $v$ is the second component.
Also the two-dimensional state is written as $\mathbf{x}=(x,y)$.

The gradient of the velocity vector is,
\begin{equation}
\nabla\mathbf{v}(\mathbf{x}) =
\begin{bmatrix}
     \frac{\partial u}{\partial x} & \frac{\partial u}{\partial y} \\ 
     \frac{\partial v}{\partial x} & \frac{\partial v}{\partial y}  
\end{bmatrix},
     \label{gradient}
\end{equation}
and  the Eulerian rate-of-strain tensor \eqref{eq:EST} is explicitly,
\begin{equation}
\mathbf{S}(\mathbf{x}) =
\begin{bmatrix}
     \frac{\partial u}{\partial x} & \tfrac{1}{2}\left(\frac{\partial u}{\partial y} + \frac{\partial v}{\partial x} \right) \\ 
     \tfrac{1}{2}\left(\frac{\partial u}{\partial y} + \frac{\partial v}{\partial x} \right) & \frac{\partial v}{\partial y}  
\end{bmatrix}.
     \label{S_matrix}
\end{equation}
The attraction rate, $s_1(\mathbf{x})$, which is the minimum eigenvalue of $\mathbf{S}(\mathbf{x})$ at the location $(\mathbf{x})$ is given analytically by,
\begin{equation}
s_1(\mathbf{x}) =
\frac{1}{2}      \left( \frac{\partial u}{\partial x} + \frac{\partial v}{\partial y} \right) - \frac{1}{2}\sqrt{\left( \frac{\partial u}{\partial x} - \frac{\partial v}{\partial y} \right)^{\hspace{-1mm}2}  
+                 \left( \frac{\partial u}{\partial y} + \frac{\partial v}{\partial x} \right)^{\hspace{-1mm}2}  },
\label{s1_analytical}
\end{equation}
where the dependence of $u$ and $v$ on $\mathbf{x}$ is understood. 
Similarly, the repulsion rate is given analytically by,
\begin{equation}
s_2(\mathbf{x}) =
\frac{1}{2}      \left( \frac{\partial u}{\partial x} + \frac{\partial v}{\partial y} \right) 
+ \frac{1}{2}\sqrt{\left( \frac{\partial u}{\partial x} - \frac{\partial v}{\partial y} \right)^{\hspace{-1mm}2}  
+                 \left( \frac{\partial u}{\partial y} + \frac{\partial v}{\partial x} \right)^{\hspace{-1mm}2}  },
\label{s2_analytical}
\end{equation}
~

\subsection{Finite-Time Lyapunov Exponent}

Based on the ODE,  
\eqref{vield_auton}, we can 
calculate,
the flow map, $\mathbf{x}_{0} \mapsto \mathbf{x}_{t} =\mathbf{F}_t(\mathbf{x}_{0})$.
The flow map, 
$\mathbf{F}_t:U \rightarrow U$, where $U \subset M$, is given by,
\begin{equation}
\mathbf{F}_t(\mathbf{x}_{0})= \mathbf{x}_{0}+\int_{0}^{t} \mathbf{v}(\mathbf{F}_\tau(\mathbf{x}_{0}))\,d\tau,\label{eq:theo:3}
\end{equation}
and is typically given numerically \cite{shadden2005definition,brunton2010fast,rypina2011investigating,pratt2014chaotic}.
Taking the gradient of the flow map $\mathbf{F}_t(\mathbf{x}_{0})$ with respect to initial conditions $\mathbf{x}_{0}$, $\nabla\mathbf{F}_t(\mathbf{x}_{0})$, the right Cauchy-Green strain tensor over a  time interval of interest can be calculated,
\begin{equation}
\mathbf{C}_t(\mathbf{x}_{0})= \nabla\mathbf{F}_t(\mathbf{x}_{0})^{T}
\nabla\mathbf{F}_t(\mathbf{x}_{0}),
\label{eq:theo:4}
\end{equation}
which is positive-definite, giving eigenvalues  which are all positive. 
In a two-dimensional flow, the two eigenvalues are ordered $\lambda_{1}<\lambda_{2}$. 
From the maximum eigenvalue, $\lambda_2$, the FTLE \cite{shadden2005definition,lekien2007lagrangian} can be defined as,
\begin{equation}
    \sigma_t(\mathbf{x}_{0}) = \frac{1}{2|t|}\log(\lambda_{2}),
    \label{FTLE}
\end{equation}
where $t$ is the (signed) elapsed time, often referred to as the integration time or evolution time in the FTLE literature.
The FTLE measures the rate of separation of two nearby fluid parcels in a flow over the time horizon $t$.
Ridges of the $\sigma_t(\mathbf{x}_{0})$ field for $t<0$ identify  regions of the flow which are the 
most  attracting over the time interval $[t,0]$. 
Furthermore, in \cite{nolan2020finite} it was demonstrated that the $-s_1(\mathbf{x}_0)$ field (with the minus sign, 
where $s_1$ is from \eqref{s1_analytical}) is the limit of the $\sigma_t(\mathbf{x}_{0})$ field as $t\rightarrow 0$. 
Thus, the attraction rate field, $s_1(\mathbf{x}_0)$,  provides a computationally inexpensive instantaneous approximation of the main attracting curves, as it is based on a single velocity snapshot; no  trajectory integration is necessary.  

\section{
On Evolution of Observations, the Koopman Operator and its Eigenfunction PDE}



The Koopman spectral analysis has become extensively popular and relevant  lately in science and engineering \cite{budivsic2012applied, kutz2016dynamic, lan2013linearization} especially for a  data-driven perspective 
for analyzing dynamical systems.  The idea is that  even a nonlinear dynamical system can be interpreted as a much simpler linear system.  However, this reinterpretation of the nonlinear dynamical system is as a linear dynamical system in a function space, which may be infinite  dimensional.  This is perhaps not a bad trade-off, infinite dimensions for linearity, and from there, computational schemes proceed to estimate the representation in terms of finite dimensional truncation of the infinite dimensional embedding.   
Often, by finite truncation the linear dynamics is approximated  by finite dimensional embedding to a linear subspace.



We first review Koopman spectral analysis in  brief details.  
Consider,  the  autonomous differential equation in
\eqref{vield_auton}.
As above, the flow for each $t\in\mathbb{R}$ (or semi-flow for $t\geq 0$) as a function, $\mathbf{x}(t)\equiv \mathbf{F}_t(\mathbf{x}_0)$ for a trajectory starting at $\mathbf{x}(0)=\mathbf{x}_0\in M$. 
There is the dynamics of the  associated Koopman operator,   often called a composition operator,  which describes  evolution of ``observables'', meaning ``measurements'' along the flow, 
\cite{budivsic2012applied, mezic2013analysis}. Rather than classically analyzing individual trajectories in the  phase space, observations measured as functions over  the space are analyzed.
These ``observation functions'',  
\begin{equation}
g:M\rightarrow {\mathbb C},
\end{equation}
are elements of a space of observation functions ${\cal F}$.  For example,
\begin{equation}
{\cal F}=L^2(M)=\{g:\int_M |g(s)|^2 ds<\infty \},
\end{equation}
is commonly used since it is particularly convenient for numerical applications that utilize the inner product associated with the Hilbert space structure 
\cite{budivsic2012applied,kutz2016dynamic, mezic2005spectral, lan2013linearization, mezic2005spectral, mezic2019spectrum}.  
We will assume scalar observation functions, but multiple scalar observation functions can just as well be considered together, ``stacked'' as a composite vector valued observation function.  

\begin{figure}
\centering
\includegraphics[width=\textwidth]{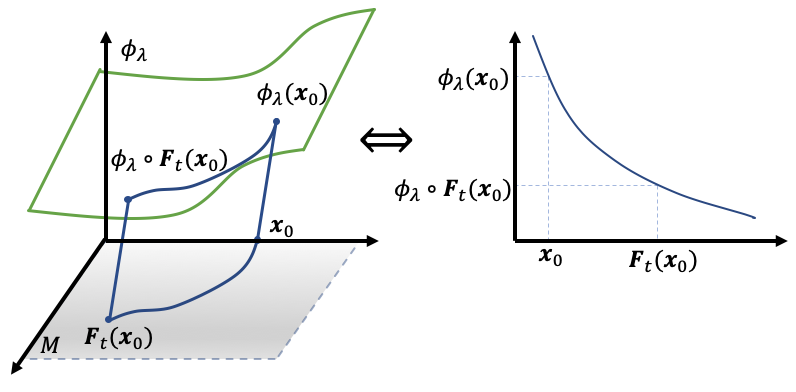}
 \caption{The action of a Koopman operator (composition operator) is to extract the value of a measurement function downstream, (Left) Eq.~(\ref{koopmandefn}), but (Right) for eigenvalue and eigenfunction pair (KEIGS) $(\lambda, \phi_\lambda(\mathbf{x}))$ that function has the property Eq.~(\ref{eq:koopman definition}) effectively interpreting the change of $\phi_\lambda$ {\it along} an orbit as if the dynamics is linear even if the flow may be nonlinear in its phase space $M$. }
\label{Fig1}
\end{figure}

The dynamics of how observation functions change over time as sampled along orbits is what the Koopman operator defines.  See Figure \ref{Fig1}.
The { Koopman Operator}, (Composition Operator), \cite{koopman1931hamiltonian, mezic2013analysis, cvitanovic2005chaos,gaspard2005chaos}, 
$\mathcal{K}_t$, associated with $\mathbf{F}_t$, is
 a (semi-) flow, stated as the following composition,
\begin{eqnarray}\label{koopmandefn}
\mathcal{K}_t:{\cal F}&\to& {\cal F}, \nonumber \\
g &\mapsto & \mathcal{K}_t[g](\mathbf{x})=g\circ \mathbf{F}_t,
\end{eqnarray}
 on the function space ${\cal F}$, for each $t\in\mathbb{R}$ (or as a semi-flow if the relation only holds for $t\geq 0$).  
That is, for each $\mathbf{x}$, we observe the value of an observable $g$ not at $\mathbf{x}$, but ``downstream" by time $t$, at $\mathbf{F}_t(\mathbf{x})$.  
See Figure \ref{Fig1}. 
Notice that for brevity we have suppressed the starting time $t_0$ in the  statement of the Koopman operator semi-flow notation, $\mathcal{K}_t$.
An  important feature of the Koopman operator is that it is linear operator on its domain, the function space $\mathcal{F}$, but at the cost of possibly being infinite dimensional, even though it may be associated with a flow $\mathbf{F}_t$ that evolves on a finite dimensional space, and indeed even due to a nonlinear vector field.

The spectral theory of Koopman operators \cite{gaspard2005chaos, budivsic2012applied, mezic2013analysis,mauroy2016global} concerns eigenfunctions and eigenvalues of the operator $K_t$.  Writing  an eigenvalue, eigenfunction pair of the Koopman operator as $(\lambda,\phi_\lambda(\mathbf{x}))$, must satisfy the equation,
\begin{equation}\label{eq:koopman definition}
\mathcal{K}_t[\phi_\lambda](\mathbf{x})=e^{\lambda t} \phi_\lambda(\mathbf{x}).
\end{equation}
See Fig.~\ref{Fig1}.
  For convenience, we will say ``KEIGs" when referring to a Koopman eigenvalue and eigenfunction pair,  $(\lambda,\phi_\lambda(\mathbf{x}))$. 
It may seem surprising if trying to make an analogy to the spectrum of matrices (finite rank operators), however, for each $\lambda$ not only is the eigenfunction $\phi_\lambda$ not unique, there are in fact {\it 
uncountably infinitely many functions associated with each $\lambda$} \cite{bollt2021geometric,korda2020optimal}.  This is true even while allowing only unit normalized eigenfunctions, to remove  the trivial idea that constant multiples of eigenfunctions are eigenfunctions.

A trending concept in the empirical study of dynamical systems has come to be the spectral decomposition of observables into eigenfunctions of the Koopman operator. 
Let a $D$-dimensional vector valued set of observables ${\bf g}(\mathbf{x})=[g_1(\mathbf{x}),...,g_D(\mathbf{x})]:M\rightarrow {\mathbb C}^D \in {\cal F}^D, $ be written as a linear combination of eigenfunctions,
\begin{equation}
    {\bf g}(\mathbf{x})=\sum_{j=1}^\infty \phi_{\lambda_i}(\mathbf{x}) {\bf v}_j
\end{equation}
where the vectors ${\mathbf v_j}\in {\mathbb C}^D$ are called Koopman modes.  Further, the power of this concept lies in the following expression that describes the dynamics of observations in terms that remind us of the linear Fourier analysis, but now of the Koopman modes.  That is if, ${\mathbf y}({\bf x}_0,t)\equiv{\mathbf g}\circ {\mathbf x}(t)= {\mathbf g}\circ \mathbf{F}_t({\mathbf x}_0)=K_t[{\mathbf g}]({\bf x}_0) $, then, 
\begin{equation}
    {\mathbf y}({\bf x}_0,t)=\sum_{j=1}^\infty  e^{\lambda_j t}\phi_{\lambda_j}({\mathbf x}_0){\mathbf v}_j.
\end{equation}
We considered the nonuniqueness of such a  decomposition, and furthermore the nature and even the cardinality of these eigenfunctions in \cite{bollt2021geometric}.  In this current work, this property leads to our goal to contrast such a Koopman spectral decomposition to other notions of coherence, such as the FTLE and iLE \cite{nolan2020finite}. 
So, we will highlight that other interesting global observations, such as the FTLE field, can be usefully presented as series in Koopman eigenfunctions.

A  fast growing literature exists concerning  data-driven approaches to construct eigenfunctions for an observed flow, given data as trajectories,  namely dynamic mode decomposition, (DMD), extended dynamic mode decomposition (EDMD) and variants \cite{schmid-2010, williams-2015,williams-2015b,li2017extended}.  
However an analytical description in terms of a PDE follows the infinitesimal generator of the Koopman operator.  
Corresponding to the statement of $\mathcal{K}_t$ as a semi-group of compositions, follows the action of the infinitesimal generator, \cite{mezic2019spectrum,lasota2013chaos,bollt2013applied},
\begin{equation}\label{ininitesimal}
     {\cal L}=\mathbf{v} \cdot \nabla,
   \end{equation}
 and so we recall, \cite{mezic2019spectrum, lasota2013chaos, bollt2018matching}, that a smooth exact eigenfunction of the Koopman operator of a given flow, for a given eigenvalue $\lambda\in {\mathbb C}$ must satisfy the following PDE, 
 \begin{equation}  
\label{kooppde}
  \mathbf{v}(\mathbf{x}) \cdot \nabla \phi_\lambda(\mathbf{x}) = \lambda \phi_\lambda(\mathbf{x}),
 \end{equation}
 if $M$ is compact, and $\phi_\lambda:M\rightarrow {\mathbb C}$, is in $C^1(M)$, or alternatively, if $\phi_\lambda\in C^2(M)$.
 Here, we will use solutions constructed directly from (\ref{kooppde}), to compare the Koopman perspective to the LCS and iLE perspective.
 
\begin{figure}[htbp]
\centering
\includegraphics[width=.55\textwidth]{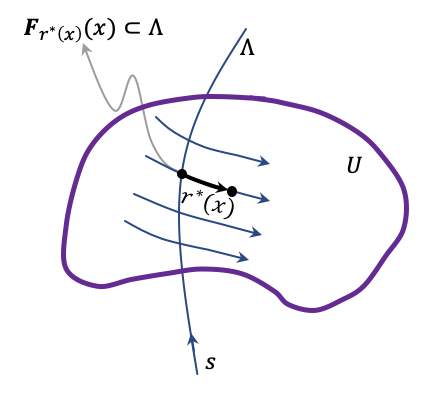}
\caption{\label{pullback} A general eigenfunction, Eq.~(\ref{gensol}) solution of Eq.~(\ref{kooppde}) is defined in terms of measuring the initial data $\protect{h(s)}$ on the data surface $\protect{\Lambda}$, and for each $\mathbf{x}$ there is a unique point where that measurement is taken.  Then the solution $\protect{\phi_\lambda}$ is a linear scaling of that measurement, by the time of the pull back.}
\end{figure}

This PDE is of a quasilinear form \cite{john1975partial}, and therefore solvable by the method of characteristics, \cite{bollt2021geometric,korda2020optimal}.  An initial data function $h:\Lambda\rightarrow {\mathbb C}$
propagates throughout an open domain in which the flow of the differential equation is defined, and respecting (\ref{kooppde}).
That an open-KEIGS pair, $(\lambda,\phi_\lambda(\mathbf{x})), \phi_\lambda:U\subset M\rightarrow {\mathbb C}$ has the form,
\begin{equation}\label{gensol}
    \phi_\lambda(\mathbf{x})=h\circ s^*(\mathbf{x})e^{\lambda r^*(\mathbf{x})},  
\end{equation}
where $r^*(\mathbf{x})$ is the ``time"-of-flight such that for a point $\mathbf{x}\in U$,  
there is an intersection in $U$ by pull back to the data surface $\Lambda$, 
\begin{equation}\label{rstar}
r^*(\mathbf{x})= \{r: \mathbf{F}_{-r}(\mathbf{x})\cap \Lambda\neq \emptyset\}.
\end{equation}  
For  each  $\mathbf{x}\in U$, 
\begin{equation}\label{sstar}
    s^*(\mathbf{x})=s\circ \mathbf{F}_{-r^*(\mathbf{x})}(\mathbf{x}),
\end{equation} is the parameterization on $\Lambda$ of that first intersection point.  This solution is valid when the orbit is non-recurrent in the open domain, during the period of interest.  See Fig.~\ref{pullback}, where we illustrate that a general solution of the eigenfunction PDE is simply a pull back along the flow through $\mathbf{x}$, to read the data on $\Lambda$, and then scale it according to the linear action of $(e^\lambda)^r$, for the ``time" it takes to pull back the point $r=r^*(\mathbf{x})$.   

\section{One-Dimensional Vector Fields}\label{1D_vector_field}

Consider  autonomous real vector fields for $x \in \mathbb{R}$,
\begin{equation}
    \dot x = f(x),
\end{equation}
The Eulerian rate-of-strain tensor is the one rate,
\begin{equation}
    s_1(x) = \frac{d f}{dx}(x)
\end{equation}
So if $s_1(x)$ is a Koopman eigenfunction, it satisfies \eqref{kooppde}, which is,
\begin{equation}
    \frac{d^2 f}{dx^2}(x) f(x) = \lambda \frac{d f}{dx}(x),
\end{equation}
for some $\lambda \in \mathbb{C}$ and $f$. 
Let $g(x)=s_1(x)$ be a real scalar function, then we can state the condition as,
\begin{equation}
    \frac{dg(x)}{dx} \int g(x)dx = \lambda g(x).
\end{equation}
That is, find a function $g(x)$ such that its derivative times its integral is proportional to the function itself. We claim that the only real function $g(x)$ that satisfies this for non-zero $\lambda$ is the zero function, $g(x)=0$. 

Therefore, for 1-dimensional flows, there is no non-trivial
instantaneous Lyapunov exponent field which is a Koopman eigenfunction.

\section{Two-Dimensional Nonlinear Saddle Flow}\label{secnonlinearsaddle}

While in 1 dimension, the iLE field is not a Koopman eigenfunction, perhaps the situation is different in 2 dimensions. As an initial 2-dimensional vector field to motivate our study, consider the following nonlinear saddle flow with cubic term,
\begin{equation}
\begin{split}
\dot x &= x, \\
\dot y &= -y + y^3,
\label{nonlinear_saddle}
\end{split}
\end{equation}
in the domain $M=\{(x,y) \in \mathbb{R}^2~\big| ~ |y|<1 \}$.
So $\mathbf{v}(\mathbf{x})=(x,-y + y^3)^T$.
These two uncoupled ordinary differential equations admit the explicit solutions,
\begin{equation}
\begin{split}
x(t) &= x e^t, \\
y(t) &= \frac{y}{\sqrt{(1 - y^2)e^{2t} + y^2}}, 
\label{nonlinear_saddle_solution}
\end{split}
\end{equation}
where the initial condition at time $0$ is $\mathbf{x}=(x,y)$.  
The right Cauchy-Green deformation tensor for a  backward integration time $t<0$, is,
\begin{equation}
\mathbf{C}_{t}(\mathbf{x}) =
\begin{bmatrix}
     e^{2t} & 0 \\ 
     0 & \frac{e^{4t}}{( (1 - y^2)e^{2t} + y^2)^3}  
\end{bmatrix},
     \label{nonlinear_saddle_C}
\end{equation}
which yields a backward-time FTLE ($t<0$), from \eqref{FTLE}, of,
\begin{equation}
    \sigma_t(\mathbf{x}) = - \frac{1}{2t} \log \Bigg( \frac{e^{4t}}{( (1 - y^2)e^{2t} + y^2)^3} \Bigg).
    \label{nonlinear_saddle_FTLE}
\end{equation}
Using Taylor series approximations for small $|t|$, the backward FTLE can be written as an expansion in $t$ for small $|t|$,
\begin{equation}
    \sigma_t(\mathbf{x}) = (1 - 3y^2) + \mathcal{O}(t).
    \label{nonlinear_saddle_FTLE_approximation_main}
\end{equation}

The backward-time FTLE can be approximated to leading order,
$\mathcal{O}(1)$, 
by the negative of the minimum eigenvalue of $\mathbf{S}(\mathbf{x})$.
The matrix $\mathbf{S}(\mathbf{x})$ is,
\begin{equation}
\mathbf{S}(\mathbf{x}) =
\tfrac{1}{2} \left( \nabla \mathbf{v} + \nabla \mathbf{v}^T  \right) 
=
\begin{bmatrix}
     1 & 0 \\ 
     0 & (- 1 + 3y^2)   
\end{bmatrix}
\end{equation}
The minimum eigenvalue is the iLE, the attracting rate \eqref{s1_analytical}, $s_{1}(\mathbf{x}) = - 1 + 3y^2$.

Using $g(\mathbf{x})=s_{1}(\mathbf{x})$ as our candidate function, we find,
\begin{equation}
\mathbf{v}(\mathbf{x}) \cdot \nabla g(\mathbf{x}) 
= 6 y^2 g(\mathbf{x}) - 12y^4
\end{equation}
which is not in the form $\lambda g(\mathbf{x})$ (as in \eqref{kooppde}), so $s_1(\mathbf{x})$ is not  a Koopman eigenfunction.
But it can 
be written 
as a {\it sum} of Koopman eigenfunctions.

Following the prescription given above 
for constructing Koopman eigenfunction using the explicit solution  \eqref{nonlinear_saddle_solution} to the ODE, we find that the Koopman eigenfunctions are of the form,
\begin{equation}
\phi_{\lambda}(\mathbf{x}) = h \left( x \sqrt{\frac{3y^2}{1-y^2}} \right) \left(\frac{3y^2}{1-y^2} \right)^{\hspace{-1mm}-\lambda/2},
\label{nonlinear_saddle_Koopman_eigenfunction}
\end{equation}
where $h : \mathbb{R} \rightarrow \mathbb{R}$, is any scalar function of $s=x \sqrt{\tfrac{3y^2}{1-y^2}}$ and $\lambda \in \mathbb{C}$ is a constant. 
One can verify directly that a function of the form \eqref{nonlinear_saddle_Koopman_eigenfunction} satisfies \eqref{kooppde}, the infinitesimal form of the eigenvalue equation.
Details are in  Appendix \ref{app:verification}.

We can now construct the scalar function $s_{1}(\mathbf{x}) = - 1 + 3y^2$ in terms of Koopman eigenfunctions.
Note that $-1$ is a Koopman eigenfunction, using $\lambda=0$ and $h(s)=-1$, a constant function, in \eqref{nonlinear_saddle_Koopman_eigenfunction}.

Also note that, via Taylor series expansion, we have,
\begin{equation}
\frac{1}{1-y^2}=1 + y^2 + y^4 + y^6 + \mathcal{O}(y^8),
\end{equation}
where $y^2 < 1$ due to our  domain $U$,
therefore,
\begin{equation}
\begin{split}
\left(\frac{3y^2}{1-y^2} \right)^{\hspace{-1mm}-\frac{\lambda}{2}} 
= (3y^2)^{-\frac{\lambda}{2}} & \big(
1 + (-\tfrac{\lambda}{2} )y^2  
+ (-\tfrac{\lambda}{2} + \tfrac{\lambda}{4} (\tfrac{\lambda}{2} +1)) y^4 \\
&+ (-\tfrac{\lambda}{2} + \tfrac{\lambda}{4} (\tfrac{\lambda}{2} +1) -\tfrac{\lambda}{12} (\tfrac{\lambda}{2} +1)(\tfrac{\lambda}{2} +2)) y^6
+  \mathcal{O}(y^8)) 
\big),
\end{split}
\end{equation}
Note that the following is a Koopman eigenfunction, using $\lambda=-2$ and $h(s)=1$,
\begin{equation}
\phi_{-2}(\mathbf{x}) = 3y^2 + 3y^4 + 3y^6 + 3y^8 + \mathcal{O}(y^{10})
\end{equation}
We can remove the $3y^4$ term by adding another  Koopman eigenfunction, using $\lambda=-4$ and $h(s)=-\tfrac{1}{3}$,
\begin{equation}
\begin{split}
\phi_{-4}(\mathbf{x})  &= -\tfrac{1}{3}
(3y^2)^2( 1 + 2y^2 + (2 + -1(-2+1))y^4 +  \mathcal{O}(y^6) ), \\
&= -3 y^4 -6y^6 -9y^8 + \mathcal{O}(y^{10}) 
\end{split}
\end{equation}
Following in a similar manner, we can remove the leading order remainder term from $\phi_{-2}(\mathbf{x}) + \phi_{-4}(\mathbf{x})$, which is $-3y^6$,  via
\begin{equation}
\begin{split}
\phi_{-6}(\mathbf{x})  
&= h(s) (3y^2)^6( 1 + 3y^2 + \mathcal{O}(y^4) ), \\
&= h(s) (27y^6 + 81y^8 + \mathcal{O}(y^{10}) ).
\end{split}
\end{equation}
Therefore, a choice of $h(s)=\tfrac{1}{9}$ will cancel out $-3y^6$.  
Following in a similar manner, we can remove the leading order remainder term from $\phi_{-2}(\mathbf{x}) + \phi_{-4}(\mathbf{x}) +  \phi_{-6}(\mathbf{x})$, which is $3y^8$, and  subsequently all higher order terms, since the leading order term of 
$\phi_{-2k}(\mathbf{x})$ is of order $y^{2k}$.  Thus we can write the term $3y^2$ as the sum of an infinite series of Koopman eigenfunctions,
\begin{equation}
3y^2 = \sum_{k=1}^{\infty} \phi_{-2k}(\mathbf{x}),
\end{equation}
where,
\begin{equation}
\phi_{-2k}(\mathbf{x}) = (-\tfrac{1}{3})^{2(k-1)} \left(\frac{3y^2}{1-y^2} \right)^{\hspace{-1.5mm}-\frac{k}{2}},
\end{equation}
for integer $k \ge 1$.
Defining  $\phi_{0}(\mathbf{x})$ as $-1$, $s_1(\mathbf{x})$, the instantaneous attracting rate, can be written exclusively in terms of these Koopman eigenfunctions,
\begin{equation}\label{examp}
s_1(\mathbf{x}) = \sum_{k=0}^{\infty} \phi_{-2k}(\mathbf{x}).
\end{equation}



We have shown that while certain special case vector fields have the property that the dominant eigenfunction is also the infinitesimal FTLEs, that is the instantaneous Lyapunov exponent (iLE), this is not the general scenario.  However, when this occurs, this would be a strong relationship between the geometric theory and the theory of evolution operators.   On the other hand, in a general scenario for a given vector field, if the corresponding eigenfunctions are dense in a space of a functions that includes the iLEs, then clearly the iLE can be written as a superposition of Koopman eigenfunctions, such as described by the
example of (\ref{examp}). Thus the geometric theory can still be interpreted in terms of the spectral theory.    In the intermediate scenario, homogeneous polynomials, (see Appendix \ref{app:arbitrary}) offer an enticing class of problems with a general relationship.  In some vector fields then, the ILE may be a finite sum of Koopman eigenfunctions.

\section{General Two-Dimensional Vector Fields}\label{2Dgeneral}

Consider a general two-dimensional autonomous vector field of the form,
\begin{equation}
\begin{split}
\dot x &= u(x,y), \\
\dot y &= v(x,y).
\label{general_2D_flow}
\end{split}
\end{equation}
where the right-hand-side functions $u$ and $v$ are as smooth as necessary in their arguments.
The gradient matrix is given by \eqref{gradient} and the attraction and repulsion rates given by \eqref{s1_analytical} and \eqref{s2_analytical}, respectively. 
We will focus on the repulsion rate, $s_2(x,y)$, but similar arguments can be applied equally to the attraction rate. 

We will proceed by seeking the conditions under which $s_2(x,y)$ is also a Koopman eigenfunction for some non-trivial eigenvalue $\lambda$, which is,
\begin{equation}
u(x,y)\pa{s_2}{x}(x,y) + v(x,y)\pa{s_2}{y}(x,y) = \lambda s_2(x,y).
\label{s2_KEIG}
\end{equation}

To simplify the calculation, we will consider only vector fields such that,
\begin{equation}
    \pa{u}{y} = -\pa{v}{x},
    \label{simplification}
\end{equation}
which, in fluid terms, would mean the vector field has zero shear component of the strain rate \cite{nolan2020finite},
in which case the repulsion rate, $s_2$, simplifies to,
\begin{equation}
s_2 = \pa{u}{x}.
\label{s2_analytical_simple}
\end{equation}
with partial derivatives,
\begin{equation}
\pa{s_2}{x} = \frac{\partial^2 u}{\partial x^2}, \quad 
\pa{s_2}{y} = \frac{\partial^2 u}{\partial x \partial y},
\end{equation}
and so the condition for the repulsion rate to be a KEIG, \eqref{s2_KEIG}, is the following condition on the vector field functions $u$ and $v$,
\begin{equation}
u\frac{\partial^2 u}{\partial x^2} + v\frac{\partial^2 u}{\partial x \partial y} = \lambda \pa{u}{x},
\label{s2_KEIG_simple}
\end{equation}
where the $(x,y)$ dependence of the vector field components is understood.
Furthermore, under the assumption \eqref{simplification}, 
the attraction rate, $s_1$, also simplifies to,
\begin{equation}
s_1= \pa{v}{y}.
\label{s1_analytical_simple}
\end{equation}
with partial derivatives,
\begin{equation}
\pa{s_1}{x} = \frac{\partial^2 v}{\partial x \partial y}, \quad 
\pa{s_1}{y} = \frac{\partial^2 v}{\partial y^2},
\end{equation}
and so the condition for the attraction rate to be a KEIG is the following condition on the vector field functions $u$ and $v$,
\begin{equation}
u\frac{\partial^2 v}{\partial x \partial y  } + v\frac{\partial^2 v}{           \partial y^2} = \lambda \pa{v}{y}.
\label{s1_KEIG_simple}
\end{equation}

\section{Polynomial Vector Fields}

It is not immediately obvious if \eqref{s2_KEIG_simple}
or \eqref{s1_KEIG_simple} admit any solutions.
We  seek to construct solutions, if they exist. 
To simplify this search, we will look initially at a special class of vector fields, polynomial vector fields.

Adopting the common notation of the polynomial vector field literature (see, e.g., \cite{schlomiuk1993algebraic}), we 
let $u(x,y)=P(x,y)$ and $v(x,y)=Q(x,y)$ be polynomials of degree $m$. 
We consider a system of the form 
\eqref{general_2D_flow} with an equilibrium point which we may place at the origin, 
\begin{equation}
\begin{split}
\dot x &= P_1(x,y) + P_2(x,y) + \cdots + P_m(x,y), \\
\dot y &= Q_1(x,y) + Q_2(x,y) + \cdots + Q_m(x,y),
\label{polynomial_2D_flow}
\end{split}
\end{equation}
where $P_k(x,y)$ and $Q_k(x,y)$ are both members of $R_k$, the homogeneous polynomials of degree $k$ (understood to be in the 2 variables $x$ and $y$).
The dimension of the vector space $R_k$ is equal to the number of monomials of degree $k$ in $x$ and $y$, and is,
${\rm dim}(R_k) = k+1$.
The number of terms in \eqref{polynomial_2D_flow} is therefore,
\begin{equation}
    2\sum_{k=1}^m {\rm dim}(R_k) = m(m+3). 
\end{equation}
We will denote the polynomials in terms of monomials,
\begin{equation}
    P_i = \sum_{j=0}^i a_{ij} x^{i-j} y^j,
    \quad
    Q_i = \sum_{j=0}^i b_{ij} x^{i-j} y^j,   
\end{equation}
where $a_{ij}$ and $b_{ij}$ are real scalar coefficients.
For \eqref{polynomial_2D_flow}, the coefficients form a space of dimension $m(m+3)$.
Based on the nonlinear saddle flow example, we expect the vector fields which satisfy \eqref{s2_KEIG_simple} or \eqref{s1_KEIG_simple}  will be some subspace of dimension $d<m(m+3)$.

\subsection{Quadratic vector fields}
Consider \eqref{polynomial_2D_flow} with $m=2$,
\begin{equation}
\begin{split}
P(x,y) &= a_{10} x + a_{11} y + a_{20} x^2 + a_{21} xy + a_{22} y^2, \\
Q(x,y) &= b_{10} x + b_{11} y + b_{20} x^2 + b_{21} xy + b_{22} y^2.
\label{quadratic_vf}
\end{split}
\end{equation}

Our assumption \eqref{simplification} implies,
\begin{equation}
\begin{split}
b_{10}  &= -a_{11},\\
b_{20}  &= -\tfrac{1}{2} a_{21},\\
b_{21}  &= -2a_{22}.
\end{split}
\end{equation} 

Therefore, the vector fields of interest are quadratic vector fields of the form,
\begin{equation}
\begin{split}
P(x,y) &= ~~a_{10} x + a_{11} y +~ ~a_{20} x^2 + ~a_{21} xy + a_{22} y^2, \\
Q(x,y) &=  -a_{11} x + b_{11} y - \tfrac{1}{2} a_{21} x^2 -2a_{22} xy + b_{22} y^2.
\label{quadratic_vf_simp}
\end{split}
\end{equation}
which is a 7-dimensional subspace of the original 10-dimensional space of quadratic vector fields. Furthermore, for this vector field, we have,
\begin{equation}
    s_2 = \pa{P}{x} = a_{10} + 2a_{20} x + a_{21} y,
    \label{s2_particular}
\end{equation}
and,
\begin{equation}
\pa{s_2}{x} = \frac{\partial^2 P}{\partial x^2} =2 a_{20}, \quad 
\pa{s_2}{y} = \frac{\partial^2 P}{\partial x \partial y} = a_{21}.
\end{equation}
Since $ s_2$ has no second order terms, the left-hand-side of \eqref{s2_KEIG_simple} must have all coefficients of second order terms identically zero, i.e.,
\begin{equation}
    \underbrace{(2a_{20}^2 - \tfrac{1}{2}a_{21}^2)}_{\rm must~be~0}x^2
    +
    \underbrace{(2a_{20}a_{21} - 2a_{21}a_{22})}_{\rm must~be~0}xy
    +
    \underbrace{(2a_{20}a_{22} + a_{21}b_{22})}_{\rm must~be~0}y^2 = 0
\end{equation}
This leads to three more algebraic conditions. Among them we choose,
\begin{equation}
\begin{split}
a_{21} &=  ~2a_{20},\\
a_{22} &= ~~a_{20},\\
b_{22} &=  -a_{20},
\end{split}
\end{equation} 
which reduces the space of possible vector fields to 4-dimensional.

The constant and linear terms of \eqref{s2_KEIG_simple} provide the following algebraic conditions,
\begin{equation}
\begin{split}
0                   &= \lambda a_{10},\\
2a_{20} (a_{10}-a_{11})x   &= \lambda 2 a_{20} x,\\
2a_{20} (a_{11} + b_{11})y   &= \lambda 2 a_{20} y.
\end{split}
\end{equation} 
This leads to conditions on $a_{10}$ and $a_{11}$,
\begin{equation}
\begin{split}
    a_{10} &= 0, \\
    a_{11} &= -\lambda, \\
    b_{11} &= 2\lambda,
\end{split}
\end{equation}
which further reduces the space of vector fields to only a 2-dimensional subspace with free parameters $(\lambda,a_{20})$.  

By construction, the repulsion rate is now a Koopman eigenfunction,
\begin{equation}
    s_2(x,y) = 2a_{20} (x + y),
    \label{s2_particular_1}
\end{equation}
with eigenvalue $\lambda$, where both $\lambda$ and $a_{20}$ are free parameters.
The quadratic vector fields which this Koopman eigenfunction corresponds to is
\begin{equation}
\begin{split}
P(x,y) &= - \lambda y + a_{20} (x+y)^2, \\
Q(x,y) &= ~~\lambda x + 2 \lambda y - a_{20} (x+y)^2.
\label{quadratic_flow_iLEKEIG_PQ}
\end{split}
\end{equation}
where $(\lambda,a_{20})\in\mathbb{R}^2$ are free parameters. 
Notice that the linear part of the vector field is skew-symmetric and that the quadratic terms in $P$ and $Q$ are opposite sign.


We note that while the repulsion rate $s_2=\tfrac{\partial P}{\partial x}$ is a KEIG, 
there is no guarantee that the attraction rate
$s_1=\tfrac{\partial Q}{\partial y}$ will be. 



\subsection{Cubic vector fields} 


For polynomial vector fields of order $m$, the $s_1$ or $s_2$ fields will be polynomials of order $m-1$, since they are based on gradients of the vector field. Thus quadratic vector fields produced a linear $s_2$ above. 
By the definition of a instantaneous Lyapunov exponent structure (iLES) \cite{nolan2020finite}, we need a ridge of $s_2$ for a repelling iLES or a trench of $s_2$ for an attracting iLES. To have this, we need an $s_1$ or $s_2$ which is quadratic. Thus, we must consider polynomial vector fields of at least cubic order.

We perform the same procedure using the same assumption \eqref{simplification}, but now for the attracting rate, $s_1$, and for vector fields satisfying \eqref{s1_KEIG_simple} which are homogeneous polynomials in $x$ and $y$ through order  $m=3$.
We augment polynomial fields of the form \eqref{polynomial_2D_flow} by allowing constant terms, $P_0=a_{00}$ and $Q_0=b_{00}$, which leads to $m(m+3)+2$ terms.
Out of the 20-dimensional space of planar cubic vector fields, we follow an approach analogous to the quadratic example above, and obtain,
\begin{equation}
\begin{split}
P(x,y) &=a_{00} +a_{10} x + a_{11} y + a_{20} (x+y)^2 + a_{20} k (x+y)^3, \\
Q(x,y) &=b_{00} -a_{11} x + b_{11} y - a_{20} (x+y)^2 - a_{20} k (x+y)^3.
\label{cubic_flow_iLEKEIG_PQ}
\end{split}
\end{equation}
where, assuming $k\ne0$, 
\begin{equation}
\begin{split}
a_{00} &= \tfrac{1}{6k}\left(a_{10}-\tfrac{a_{20}}{3k}\right) - b_{00}, \\
a_{11} &= \tfrac{1}{2} \left(a_{10}+\tfrac{a_{20}}{3k}\right), \\
b_{11} &=                       -\tfrac{a_{20}}{3k}, 
\label{cubic_flow_iLEKEIG_params}
\end{split}
\end{equation}
where $(a_{10},k,a_{20},b_{00})\in\mathbb{R}^4$ are free parameters. Notice that the linear part of the vector field is skew-symmetric and that the quadratic and cubic terms in $P$ and $Q$ are opposite sign.

For this cubic vector field, the attracting rate function,
\begin{equation}
    s_1 = \tfrac{\partial Q}{\partial y} = -\tfrac{a_{20}}{3k} -2a_{20} (x+y) -3a_{20} k(x+y)^2,
\end{equation}
is a Koopman eigenfunction with eigenvalue,
\begin{equation}
    \lambda = a_{10} - \tfrac{a_{20}}{3k}.
\end{equation}
Because the attracting rate field is quadratic, it can potentially have a ridge. 
We re-write the vector field in terms of $\lambda$ and $c=-2ka_{20}$,
\begin{equation}
\begin{split}
P(x,y) &=a_{00} + ~~(\lambda - \tfrac{1}{6}\tfrac{c}{k^2}) x + (\tfrac{1}{2}\lambda - \tfrac{1}{6}\tfrac{c}{k^2}) y - \tfrac{1}{2}\tfrac{c}{k} (x+y)^2 - \tfrac{1}{2}c (x+y)^3, \\
Q(x,y) &=(\tfrac{1}{6}\tfrac{\lambda}{k}-a_{00}) -(\tfrac{1}{2}\lambda - \tfrac{1}{6}\tfrac{c}{k^2}) x + \tfrac{1}{6}\tfrac{c}{k^2} y + \tfrac{1}{2}\tfrac{c}{k} (x+y)^2 + \tfrac{1}{2}c (x+y)^3,
\label{cubic_flow_iLEKEIG}
\end{split}
\end{equation}
where $(\lambda,c,k,a_{00})\in\mathbb{R}^4$ are free parameters. 
The attracting rate, 
\begin{equation}
    s_1 = \tfrac{1}{6}\tfrac{c}{k^2} +\tfrac{c}{k} (x+y) +\tfrac{3}{2} c (x+y)^2,
\end{equation}
is a Koopman eigenfunction with eigenvalue $\lambda$ for vector field \eqref{cubic_flow_iLEKEIG}.

\subsection{Cubic vector field example}

As an explicit example of \eqref{cubic_flow_iLEKEIG}, consider the following, where 
$(\lambda,c,k,a_{00})=(2,\tfrac{2}{3},-\tfrac{1}{3},-2)$,
\begin{equation}
\begin{split}
\dot x &= -2 +  x + (x+y)^2 - \tfrac{1}{3} (x+y)^3, \\
\dot y &=~~1 +  y - (x+y)^2 + \tfrac{1}{3} (x+y)^3.
\label{cubic_flow_iLEKEIG_example}
\end{split}
\end{equation}
The attracting rate for  this vector field,
\begin{equation}
    s_1 =  1 - 2 (x+y) +(x+y)^2 = \left( (x+y) - 1\right)^2,
\end{equation}
is a Koopman eigenfunction with eigenvalue $\lambda=2$. We note that if we perform a linear  transformation to new variables,
\begin{equation}
    \begin{bmatrix}
     r \\ 
     s   
    \end{bmatrix} = T 
        \begin{bmatrix}
     x \\ 
     y   
    \end{bmatrix}
\quad {\rm where} \quad
T =
\begin{bmatrix}
      ~~1 & 1\\ 
       -1 & 1   
    \end{bmatrix},
\end{equation}
we get a vector field that does not have the same symmetries,
\begin{equation}
\begin{split}
\dot r &= -1 + r, \\
\dot s &=~~3 + s- 2r^2 +\tfrac{2}{3} r^3.
\label{cubic_flow_iLEKEIG_example_alt}
\end{split}
\end{equation}
but nonetheless has the same attracting rate, 
\begin{equation}
    s_1 =  ( r-1)^2,
\end{equation}
which is  a Koopman eigenfunction of the new velocity field,
with the same eigenvalue $\lambda =2$.

\subsection{Cubic vector field transformation to simplify}

For general parameters, the same linear transformation $T$ leads to the vector field \eqref{cubic_flow_iLEKEIG} in the $(r,s)$ variables as follows,
\begin{equation}
\begin{split}
\dot r &= \tfrac{1}{6}\tfrac{\lambda}{k} + \tfrac{1}{2}\lambda r, \\
\dot s &=\left( \tfrac{1}{6}\tfrac{\lambda}{k}-2a_{00} \right) + 
\left( - \lambda + \tfrac{1}{3}\tfrac{c}{k^2} \right) r +
\tfrac{1}{2}\lambda s +\tfrac{c}{k}r^2 +c r^3.
\label{cubic_flow_iLEKEIG_example_alt2}
\end{split}
\end{equation}
As the $\dot r$ equation is uncoupled from $s$, it can be solved in closed form. For an initial condition $r_0$ at time $t_0$, we have,
\begin{equation}
    r(t) =  \tfrac{1}{a} \left[ (\bar a + a r_0)e^{a(t-t_0)}-\bar a \right]
\end{equation}
where $a=\tfrac{1}{2}\lambda$ and $\bar a=\tfrac{1}{6}\tfrac{\lambda}{k}$. 
Therefore the $s$ ODE can be considered a time-dependent ODE in 1D.

Notice that $\bar r =-\bar a/a=-\tfrac{1}{3k}$ is an invariant manifold, i.e., $\dot r =0$ (and this is also the location of the attracting rate iLE ridge for $k=\pm \tfrac{1}{3}$).

It turns out we can transform the nonlinear ODE \eqref{cubic_flow_iLEKEIG_example_alt2} into a form
where it can be solved analytically. 
We  translate the vector field so the new origin is the single equilibrium point of \eqref{cubic_flow_iLEKEIG_example_alt2},
\begin{equation}
\begin{split}
\bar r &= - \tfrac{1}{3k}, \\
\bar s &= - \tfrac{2}{\lambda} \left(  \tfrac{1}{2}\tfrac{\lambda}{k}  - \tfrac{1}{27}\tfrac{c}{k^3} -2 a_{00} \right),
\label{cubic_flow_iLEKEIG_example_alt2_eq}
\end{split}
\end{equation}
by defining,
\begin{equation}
\begin{split}
x_1 &= r - \bar r, \\
x_2 &= s - \bar s,
\label{cubic_flow_iLEKEIG_example_alt2_translate}
\end{split}
\end{equation}
where the vector field in $(x_1,x_2)$ is,
\begin{equation}
\begin{split}
\dot x_1 &= ~~\tfrac{\lambda}{2} x_1, \\
\dot x_2 &=  -\lambda x_1 + \tfrac{\lambda}{2} x_2 + c x_1^3.
\label{cubic_flow_iLEKEIG_example_alt2_translate_ode}
\end{split}
\end{equation}
Notice this transformed system 
depends on only 2 parameters, $(\lambda,c)$.
Furthermore, the vector field \eqref{cubic_flow_iLEKEIG_example_alt2_translate_ode} does not satisfy the simplifying assumption \eqref{simplification} which held in the original variables.
We note that the origin of \eqref{cubic_flow_iLEKEIG_example_alt2_translate_ode} is the only equilibrium point and is either a stable or unstable node, depending on whether $\lambda$ is negative or positive, respectively. 

The system has an $s_1$ field that is a KEIG, $\phi_{\lambda}$, with eigenvalue $\lambda$, 
\begin{equation}
    s_1 = -\tfrac{3}{2}cx_1^2.
    \label{s1_w}
\end{equation}
which has a trench at $x_1=0$ for $c<0$, which is a necessary condition for being an attracting iLES \cite{nolan2020finite}. From \eqref{cubic_flow_iLEKEIG_example_alt2_translate_ode}, we see that the set $\{x_1=0\}$, the $x_2$-axis, is also an invariant manifold. For $\lambda<0$ $(\lambda>0)$, the $x_2$-axis is the fastest direction along the  stable (unstable) invariant manifold of the origin.


\subsection{Properties of KEIGs of 2-dimensional cubic vector fields}

We discuss some of the properties of KEIGs of cubic vector fields in two dimensions.
The nonlinear ODE \eqref{cubic_flow_iLEKEIG_example_alt2_translate_ode} is, nonetheless, in a form which admits Carleman linearization
\cite{carleman1932}.  

Augment the system with a third nonlinear variable, $x_1^3$.
Then we have a linear ODE,
\begin{equation}
    \mathbf{y} = \begin{bmatrix}
         y_1 \\ y_2 \\ y_3
    \end{bmatrix} 
    =
    \begin{bmatrix}
         x_1 \\ x_2 \\ x_2^3
    \end{bmatrix}
    \quad \Longrightarrow \quad
    \begin{bmatrix}
        \dot y_1 \\ \dot y_2 \\ \dot y_3
    \end{bmatrix}  =   
  \underbrace{\begin{bmatrix}
        \tfrac{\lambda}{2} & 0 & 0 \\
        -\lambda & \tfrac{\lambda}{2} & c \\
        0 & 0 & \tfrac{3}{2}\lambda
    \end{bmatrix}
    }_{\mathbf{A}}
  \begin{bmatrix}
         y_1 \\ y_2 \\ y_3
    \end{bmatrix}, 
\end{equation}
where the $\mathbf{y}$ space can be interpreted as a three-dimensional Koopman observable vector space \cite{brunton2019data}. 
The linear ODE in $\mathbf{y}$ admits an analytical solution,
\begin{equation}
    \mathbf{y}(t) = e^{\mathbf{A}t}\mathbf{y}(0),
\end{equation}
which implies $x_1(t)=e^{\lambda t/2}x_1(0)$. 

Using $s_1$ in \eqref{s1_w} as the observable function in \eqref{koopmandefn}, we can explicitly find that,
\begin{equation}
    \begin{split}
        s_1(t) 
        &= -\tfrac{3}{2}c x_1(t)^2, \\
        &= e^{\lambda t} \left( -\tfrac{3}{2}c x_1(0)^2 \right), \\
        & =(e^{\lambda})^t s_1(0),
    \end{split}
\end{equation}
so following individual trajectories, the $s_1$ field at time $t$ is the same as the $s_1$ field at the initial time $0$ multiplied by a factor $(e^{\lambda})^t$ as we expect from an observable that is also a KEIG; see  \eqref{eq:koopman definition}.

\begin{figure}[h!]
	\centering
	\includegraphics[width=0.8\textwidth]{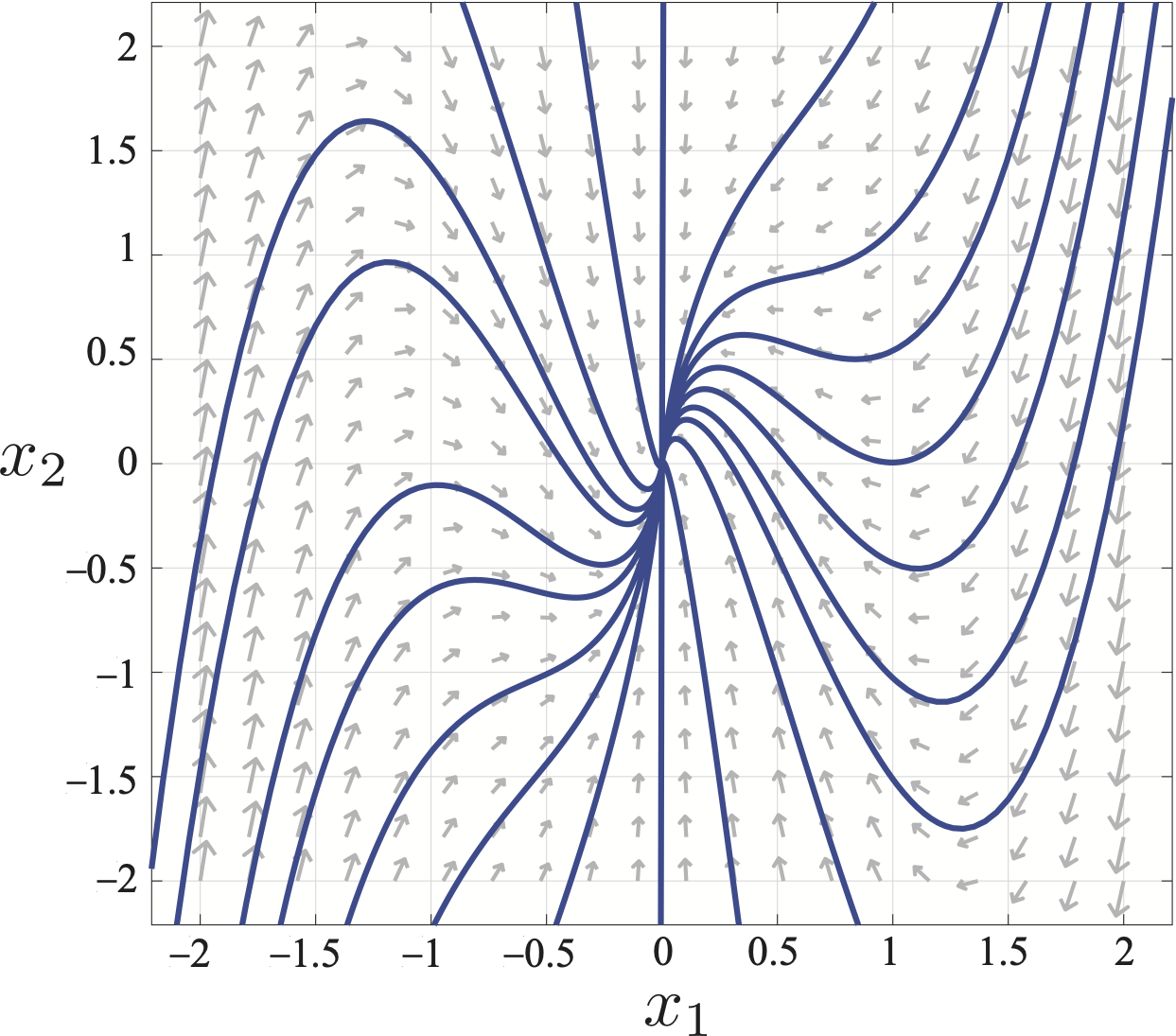}
	\caption{Phase portrait  for vector field  \eqref{cubic_flow_iLEKEIG_example_alt2_translate_ode} in $(x_1,x_2)$ coordinates for parameters $(\lambda,c)=(-1,-1)$. 
	The attracting rate field for this vector field is also a Koopman eigenfunction (KEIG). The $x_2$-axis $\{x_1=0\}$, as trench of the attracting rate field, is an attracting instantaneous Lyapunov exponent structure (iLES).}
	\label{fig:phase_portrait_lam_c_1}
\end{figure}

We show an example phase portrait for $(\lambda,c)=(-1,-1)$ in Figure \ref{fig:phase_portrait_lam_c_1}.
For these parameters, the $x_2$-axis, $\{x_1=0\}$, is an attracting iLES.

Note that $s_1$ is not unique as a Koopman eigenfunction with eigenvalue $\lambda$. 
We can verify that $x_1$ is itself a KEIG with eigenvalue $\tfrac{\lambda}{2}$, and due to the theorem of \cite{bollt2021geometric}, any power of $x_1$ is also a KEIG, 
in particular, $x_1^2$ is a KEIG with eigenvalue $\lambda$, as is $\gamma x_1^2$ for any constant $\gamma$.

\subsection{A family of polynomial vector fields with a KEIG attracting rate}

Based on the form of \eqref{cubic_flow_iLEKEIG_example_alt2_translate_ode}, we can construct a family of vector fields which all have the attracting rate as a KEIG,
\begin{equation}
\begin{split}
\dot x_1 &= ~~\tfrac{\lambda}{2} x_1, \\
\dot x_2 &=  -\lambda x_1 + \tfrac{\lambda}{2} x_2 + c x_1^3
+ c_4 x_1^4 + c_5 x_1^5 + \cdots 
\label{cubic_flow_iLEKEIG_example_alt2_translate_ode_fam}
\end{split}
\end{equation}
where we include in $\dot x_2$ only additional powers of $x_1$, with additional parameters $c_4$, $c_5$, etc. 
The vector field \eqref{cubic_flow_iLEKEIG_example_alt2_translate_ode_fam}
has the attracting rate,
\begin{equation}
    s_1 = -\tfrac{1}{2} \left( 3 c x_1^2 + 4 c_4 x_1^3 + 5 c_5 x_1^4 + \cdots \right).
    \label{s1_fam}
\end{equation}
which is a KEIG with eigenvalue $\lambda$. As before, the $x_2$-axis is (i) an attracting iLES for $c<0$ and (ii) the fastest direction along the  stable (unstable) invariant manifold of the origin for $\lambda<0$ ($\lambda >0$).

\section{Conclusion}\label{s:conclusion}

In the present work, we construct vector fields with the property that either the attracting rate or repulsion rate is a Koopman eigenfunction. 
We find that in 1 dimension, this is not possible, but in 2 dimensions, it is. 
We consider cubic 2-dimensional vector fields, and find a 2-parameter family of systems which have an attracting rate KEIG.
It turns out that these systems, while nonlinear, have the 
unusual property that Carleman linearization truncates at finite order, 
allowing us to find the exact analytical solution of the flow map using linear methods.
It is not obvious why the search for a vector field where the attracting rate is a KEIG would lead to this property.

The further investigation of vector fields which have attracting and repulsion rate fields as KEIGs is left as future work.
Returning to the title question of this paper, ``Is the finite-time Lyapunov exponent field a Koopman eigenfunction?" while the answer is yes in some special cases, in general, the answer is, no.

\section*{Acknowledgments}
S.D.R.  was supported in part by National Science Foundations grants 1821145 and 1922516,
as well as the National Aeronautics and Space Administration under Grant No.\ 80NSSC20K1532 
issued through the Interdisciplinary Research in Earth Science (IDS) and Biological Diversity programs.
E.M.B. received funding from the Army Research Office (N68164-EG), ONR and also DARPA. We thank A.J.
for a careful reading of this manuscript.

\section*{Appendices}

\begin{appendices}

\section{Verification of the form of the Koopman eigenfunctions}\label{app:verification}

Continuing from the discussion in section \ref{secnonlinearsaddle},
we follow the prescription given by Bollt \cite{bollt2021geometric} for constructing Koopman eigenfunction using the explicit solution to the ODE, \eqref{nonlinear_saddle_solution}. 
We get that Koopman eigenfunctions are of the form,
\begin{equation}
\phi_{\lambda}(\mathbf{x}) = h \left( x \sqrt{\frac{3y^2}{1-y^2}} \right) \left(\frac{3y^2}{1-y^2} \right)^{\hspace{-1.5mm}-\frac{\lambda}{2}},
\nonumber 
\end{equation}
where $h : \mathbb{R} \rightarrow \mathbb{R}$, is any scalar function, $h(s)$, of $s=x \sqrt{\tfrac{3y^2}{1-y^2}}$, and $\lambda \in \mathbb{C}$ is a constant.  One can verify directly that a function of the form \eqref{nonlinear_saddle_Koopman_eigenfunction} satisfies the infinitesimal form of the eigenvalue equation,
\begin{equation}
\mathbf{v}(\mathbf{x}) \cdot \nabla \phi_{\lambda}(\mathbf{x})  = \lambda \phi_{\lambda}(\mathbf{x}),
\nonumber 
\end{equation}
We start by taking the partial derivatives of \eqref{nonlinear_saddle_Koopman_eigenfunction}:
\begin{equation}
\tfrac{\partial}{\partial x} \phi_{\lambda}(\mathbf{x})  = 
h^{\prime} \left(\frac{3y^2}{1-y^2} \right)^{\hspace{-1.5mm}-\frac{\lambda}{2}} \sqrt{\frac{3y^2}{1-y^2} },
\label{grad_phi_x}
\end{equation}
and
\begin{equation}
\begin{split}
\tfrac{\partial}{\partial y} \phi_{\lambda}(\mathbf{x})  = 
h^{\prime} x \tfrac{1}{2} &\frac{1}{\sqrt{\frac{3y^2}{1-y^2}}} \left(\frac{3y^2}{1-y^2} \right)^{\hspace{-1.5mm}-\frac{\lambda}{2}} \frac{\partial}{\partial y}\left(\frac{3y^2}{1-y^2} \right)  \\
&+ h \tfrac{\lambda}{2} \frac{1}{\left(\frac{3y^2}{1-y^2}\right)} \left(\frac{3y^2}{1-y^2} \right)^{\hspace{-1.5mm}-\frac{\lambda}{2}}
\frac{\partial}{\partial y}\left(\frac{3y^2}{1-y^2} \right)
\label{grad_phi_y}
\end{split}
\end{equation}
Now take the dot product of the gradient, 
\[
\nabla \phi_{\lambda}(\mathbf{x}) = \left( \tfrac{\partial}{\partial x} \phi_{\lambda}(\mathbf{x}),\tfrac{\partial}{\partial y} \phi_{\lambda}(\mathbf{x}) \right),
\]
with $\mathbf{v}(\mathbf{x})$ from \eqref{nonlinear_saddle}.
\begin{equation}
\begin{split}
\mathbf{v}(\mathbf{x}) \cdot \nabla \phi_{\lambda}(\mathbf{x}) &=  
h^{\prime} x \left(\frac{3y^2}{1-y^2} \right)^{\hspace{-1.5mm}-\frac{\lambda}{2}} \sqrt{\frac{3y^2}{1-y^2} } \\
&+
h^{\prime} x \tfrac{1}{2} \frac{1}{\sqrt{\frac{3y^2}{1-y^2}}} \left(\frac{3y^2}{1-y^2} \right)^{\hspace{-1.5mm}-\frac{\lambda}{2}} \frac{\partial}{\partial y}\left(\frac{3y^2}{1-y^2} \right) (-y+y^3)\\ 
&-
h \tfrac{\lambda}{2} \frac{1}{\left(\frac{3y^2}{1-y^2}\right)} \left(\frac{3y^2}{1-y^2} \right)^{\hspace{-1.5mm}-\frac{\lambda}{2}}
\frac{\partial}{\partial y}\left(\frac{3y^2}{1-y^2} \right) (-y+y^3)
\label{v_grad_dot_product}
\end{split}
\end{equation}
Now, note that,
\begin{equation}
\begin{split}
\frac{1}{2} \frac{1}{\sqrt{\frac{3y^2}{1-y^2}}} \frac{\partial}{\partial y}\left(\frac{3y^2}{1-y^2} \right) (-y+y^3)
& = \frac{1}{2} \sqrt{\frac{1-y^2}{3y^2}} \frac{6y}{(1-y^2)^2}(-y+y^3) \\
& = \frac{\sqrt{3} (-y)}{(1-y^2)^{\frac{3}{2}}} (1-y^2) \\
& = -\sqrt{\frac{3y^2}{1-y^2} } 
\end{split}
\end{equation}
Therefore, the first two terms of \eqref{v_grad_dot_product} cancel out. The remaining term is,
\begin{equation}
\begin{split}
-h \left(\frac{3y^2}{1-y^2} \right)^{\hspace{-1.5mm}-\frac{\lambda}{2}} \frac{\lambda}{2} \frac{(1-y^2)}{3y^2} \frac{6y}{(1-y^2)^2}(-y)(1-y^2)
&=\lambda h \left(\frac{3y^2}{1-y^2} \right)^{\hspace{-1.5mm}-\frac{\lambda}{2}},\\
&=\lambda \phi_{\lambda}(\mathbf{x}),\\
\end{split}
\end{equation}
We have therefore verified that the Koopman eigenfunctions of the nonlinear saddle, \eqref{nonlinear_saddle}, are of the form \eqref{nonlinear_saddle_Koopman_eigenfunction}.

\section{Arbitrary functions $g(\mathbf{x})$ can be written as an infinite series of Koopman eigenfunctions}\label{app:arbitrary}

For the nonlinear saddle system in section \ref{secnonlinearsaddle},
let $h(s)=s^n$, then,
\begin{equation}
\begin{split}
\phi_{n,\lambda}(\mathbf{x}) 
&= h \left( x \sqrt{\tfrac{3y^2}{1-y^2}} \right) \left(\tfrac{3y^2}{1-y^2} \right)^{\hspace{-1.5mm}-\frac{\lambda}{2}},\\
&= \left( x \sqrt{\tfrac{3y^2}{1-y^2}} \right)^{\hspace{-1.5mm}n} \left(\tfrac{3y^2}{1-y^2} \right)^{\hspace{-1.5mm}-\frac{\lambda}{2}},\\
&=  x^n \left(\tfrac{3y^2}{1-y^2} \right)^{\hspace{-1.5mm}\frac{n-\lambda}{2}},\\
\end{split}
\end{equation}
So,
\begin{equation}
\begin{split}
\phi_{1,1}(\mathbf{x}) &= x,\\
\phi_{2,2}(\mathbf{x}) &= x^2,\\
\phi_{n,n}(\mathbf{x}) &= x^n,\\
\phi_{1,0}(\mathbf{x}) &= x  y \sqrt{\tfrac{3}{1-y^2} } ,\\
\phi_{2,1}(\mathbf{x}) &= x^2  y \sqrt{\tfrac{3}{1-y^2} } ,\\
\phi_{2,0}(\mathbf{x}) &= x^2  y^2 \left(\tfrac{3}{1-y^2} \right)  ,\\
\phi_{0,-1}(\mathbf{x}) &= y \sqrt{\tfrac{3}{1-y^2} }, \\
& =  \sqrt{3} y (1 + y^2 + y^4 + y^6 + \mathcal{O}(y^8) )^{\frac{1}{2}} ,\\
& =  \sqrt{3} y (1 + \tfrac{1}{2}y^2 - \tfrac{1}{8}y^4 + \tfrac{1}{16}y^6 + \mathcal{O}(y^8) ),\\
\phi_{0,-3}(\mathbf{x}) &= y^3 \left(\tfrac{3}{1-y^2} \right)^{\hspace{-1.5mm}\frac{3}{2}}, \\
& =  (3)^{\frac{3}{2}} y^3 (1 + y^2 + y^4 + y^6 + \mathcal{O}(y^8) )^{\frac{3}{2}} ,\\
& =  (3)^{\frac{3}{2}} y^3 (1 + \tfrac{3}{2}y^2 + \tfrac{3}{8}y^4 - \tfrac{1}{16}y^6 + \mathcal{O}(y^8) ),
\end{split}
\end{equation}
So $y$ is an infinite series of Koopman eigenfunctions,
\begin{equation}
\begin{split}
y &=(3)^{-\frac{1}{2}} \phi_{0,-1}(\mathbf{x}) - 2(3)^{-\frac{3}{2}}\phi_{0,-3}(\mathbf{x}) + \tfrac{25}{8}(3)^{-\frac{5}{2}}\phi_{0,-5}(\mathbf{x}) + \cdots
\end{split}
\end{equation}
In a similar way, we can write $y^n$ for an integer $n$ which is either even or odd, as an infinite series of Koopman eigenfunctions, and the same for $x^n y^m$.
So 
any function $g(\mathbf{x})$ which is a homogeneous polynomial of order $n\in\mathbb{Z}$ in $x$ and $y$ can be written in terms of an infinite series of Koopman eigenfunctions.  
This means any function which admits a Taylor series expansion can be written in terms of an infinite series of Koopman eigenfunctions.

\end{appendices}

\clearpage

\bibliography{main}

\end{document}